\input amstex
\documentstyle{amsppt}
\pagewidth{5.4in}
\pageheight{7.6in}
\magnification=1200
\TagsOnRight
\NoRunningHeads
\topmatter
\title
\bf Some results for the Perelman LYH-type inequality
\endtitle
\author
Shu-Yu Hsu
\endauthor
\affil
Department of Mathematics\\
National Chung Cheng University\\
168 University Road, Min-Hsiung\\
Chia-Yi 621, Taiwan, R.O.C.\\
e-mail:syhsu\@math.ccu.edu.tw
\endaffil
\date
May 12, 2008
\enddate
\address
e-mail address:syhsu\@math.ccu.edu.tw
\endaddress
\abstract
Let $(M,g(t))$, $0\le t\le T$, $\partial M\ne\phi$, be a compact 
$n$-dimensional manifold, $n\ge 2$, with metric $g(t)$ evolving by the Ricci 
flow such that the second fundamental form of $\partial M$ with respect to 
the unit outward normal of $\partial M$ is uniformly bounded below on 
$\partial M\times [0,T]$. We will prove a global Li-Yau gradient estimate 
for the solution of the generalized conjugate heat equation on $M\times 
[0,T]$. We will give another proof of Perelman's Li-Yau-Hamilton type 
inequality for the fundamental solution of the conjugate heat equation 
on closed manifolds without using the properties of the reduced distance. 
We will also prove various gradient estimates for the Dirichlet 
fundamental solution of the conjugate heat equation. 
\endabstract
\endtopmatter
\NoBlackBoxes
\define \pd#1#2{\frac{\partial #1}{\partial #2}}
\define \1{\partial}
\define \2{\overline}
\define \3{\varepsilon}
\define \4{\widetilde}
\define \5{\underline}
\define \ov#1#2{\overset{#1}\to{#2}}
\define \oa#1{\overset{a}\to{#1}}
\document

In \cite{P} Perelman stated a differential Li-Yau-Hamilton type inequality
for the fundamental solution of the conjugate heat equation on closed
manifolds evolving by the Ricci flow. More precisely let $M$ be a closed 
manifold with metric $g(t)$, $0\le t\le T$, evolving by the Ricci flow,
$$
\frac{\1}{\1 t} g_{ij}=-2R_{ij}\tag 0.1
$$
in $M\times (0,T)$. Let $p\in M$ and
$$
u=\frac{e^{-f}}{(4\pi\tau)^{\frac{n}{2}}}\tag 0.2
$$
be the fundamental solution of the conjugate heat equation 
$$
u_t+\Delta u-Ru=0\tag 0.3
$$
in $M\times (0,T)$ where $\tau=T-t$ and $R=R(\cdot,t)$ is the scalar curvature
of $M$ with respect to the metric $g(t)$ with 
$$
\lim_{t\nearrow T}u=\delta_p\tag 0.4
$$ 
in the distribution sense where $\delta_p$ is the
delta mass at $p$. Let
$$
v=[\tau (2\Delta f-|\nabla f|^2+R)+f-n]u\tag 0.5
$$
where $\tau=T-t$. Then
$$
v(x,t)\le 0\quad\text{ in }M\times (0,T).\tag 0.6
$$
This result was used by Perelman to give a proof of the pseudolocality
theorem in section 10 of \cite{P} which more or less said that almost Euclidean
regions of large curvature in closed manifold with metric evolving by 
Ricci flow remain localized. Perelman gave a sketch of the proof of (0.6) 
in \cite{P} and a detailed proof 
of it using properties of reduced distance was later given by L.~Nei 
\cite{N3}. This result was generalized by L.~Nei \cite{N1}, \cite{N2}, to 
the case of the linear heat equation and by A.~Chau, L.F.~Tam, and C.~Yu 
\cite{CTY} to complete manifold with uniformly bounded curvatures. 

Let $(M,g(t))$, $0\le t\le T$, $\1 M\ne\phi$, be a compact $n$-dimensional 
manifold, $n\ge 2$, with metric $g(t)$ evolving by the Ricci flow 
such that the second fundamental form II of $\1 M$ with respect to the unit 
outward normal $\1/\1\nu$ of $\1 M$ is uniformly bounded below on $\1 M\times 
[0,T]$. In this paper we will use a variation of the method of P.~Li, 
S.T.~Yau, \cite{LY} and J.~Wang \cite{W} to prove a global Li-Yau gradient 
estimate for the solution of the generalized conjugate heat equation on 
such manifold with Neumann boundary condition. 

We obtain a similar type of global gradient estimate for the solution 
of the generalized conjugate heat equation on closed manifold with metric 
evolving by the Ricci flow. As a consequence we obtain another proof of 
Perelman's Li-Yau-Hamilton type inequality for the fundamental solution 
of the conjugate heat equation on closed manifolds without using the 
properties of the reduced distance. 

We will also prove various  gradient estimates for the Dirichlet 
fundamental solution of the conjugate heat equation. Note that localized 
Li-Yau estimate 
for the conjugate heat equation on compact manifolds with metric evolving 
by the Ricci flow was also proved by S.~Kuang and Q.S.~Zhang in \cite{KZ}. 
We refer the readers to the paper \cite{H} by R.S.~Hamilton for the recent 
results on Ricci flow and the book \cite{CLN} by B.~Chow, P.~Lu and L.~Ni 
for the basics of Ricci flow.

The plan of the paper is as follows. In section 1 we will prove a global 
Li-Yau gradient estimate for the solution of the generalized conjugate heat 
equation on compact manifolds with boundary and on closed manifolds. In 
section 2 we will give another proof of Perelman's Li-Yau-Hamilton 
type inequality on closed manifolds without using the properties of reduced 
distance. In section 3 we will generalize a result of Q.S.~Zhang \cite{Z}
to local gradient estimates for the solutions of generalized conjugate heat
equation. In section 4 we will prove the gradient estimates for the 
Dirichlet fundamental solution of the conjugate heat 
equation. 

We start with some definitions. Let $\nabla^t$ and $\Delta^t$ be the 
covariant derivative and Laplacian with respect to the metric $g(t)$. 
When there is no ambiguity, we will drop the superscript and write $\nabla$,
$\Delta$, for $\nabla^t$, $\Delta^t$, respectively. For any $r>0$, 
$x_0\in M$, $0<t_1\le t_0\le T$, let $B_r(x_0)$ be the 
geodesic ball with center $x_0$ and radius $r$ with respect to the 
metric $g(0)$ and $Q_{r,t_1}(x_0,t_0)=B_r(x_0)\times [t_0-t_1,t_0]$.
Let $dV_t$ be the volume element with respect to the metric $g(t)$
and let $V_x^t(r)=\text{Vol}_{g(t)}(B_r(x))$, $V_x(r)=V_x^0(r)$.
For any $x_1,x_2\in M$, let $r(x_1,x_2)$ be the distance between $x_1$ 
and $x_2$ with respect to the metric $g(0)$.

We also recall a definition of R.~Chen \cite{C}. For any $x\in M$ let $r(x)$ 
be the distance of $x$ from $\1 M$ with respect to $g(0)$. We say that $\1 M$ 
satisfies the interior rolling $R$-ball condition if for for any $p\in\1 M$,
there exists a geodesic ball $B_{R/2}(q)\subset M$ with center at $q\in M$ 
and radius $R/2$ respect to the metric $g(0)$ such that $\{p\}
=\2{B_{R/2}(q)}\cap\1 M$.

$$
\text{Section 1}
$$

In this section unless stated otherwise, we will let $(M,g(t))$, 
$0\le t\le T$, $\1 M\ne\phi$, be a compact $n$-dimensional manifold, 
$n\ge 2$, with metric $g(t)$ satisfying
$$
\frac{\1}{\1 t}g_{ij}=2h_{ij}\quad\text{ on }M\times [0,T]\tag 1.1
$$ 
where $h_{ij}(x,t)$ is a smooth family of symmetric tensors on $M$. We will 
assume that the second fundamental form II of $\1 M$ with respect to the unit 
outward normal $\1/\1\nu$ of $\1 M$ and metric $g(t)$ is uniformly bounded 
below by $-H$ for all $0\le t\le T$ and 
$$
|Rm|\le k_0\quad\text{ on }M\times [0,T]\tag 1.2
$$
for some constants $H>0$ and $k_0>0$. Let $u$ be a positive solution of 
$$\left\{\aligned
&u_t=\Delta^tu-qu\quad\text{ in }M\times [0,T]\\
&\frac{\1 u}{\1\nu}=0\qquad\qquad\text{ on }\1 M\times (0,T)
\endaligned\right.\tag 1.3
$$
where $q(x,t)$ is a smooth function of $M\times [0,T]$.

In this section we will prove a global Li-Yau gradient estimate for the 
solution of (1.3) on $M\times (0,T)$. We start with an algebraic lemma.

\proclaim{\bf Lemma 1.1}
Let $A, B\in\Bbb{R}$, $A\ge 0$, be constants satisfying $B\le A/\alpha$ for 
some constant $\alpha>1$. For any $0<\rho<1$, let $I(\rho)=(A-B)^2-\rho A^2$. 
Then there exists a constant $\rho\in (0,1)$ such that
$$
I(\rho)\ge\frac{1}{\alpha^2}(A-\alpha B)^2.\tag 1.4
$$ 
\endproclaim
\demo{Proof}
We divide the proof into two cases.

\noindent $\underline{\text{\bf Case 1}}$: $B>0$.

Then $B^2\le A^2/\alpha^2$. By direct computation for any $0<\sigma<1$,
$$\align
I(\rho)=&(A-\alpha B)^2+2(\alpha-1)B(A-\alpha B)+(\alpha -1)^2B^2-\rho A^2\\
=&(1-\sigma)(A-\alpha B)^2+J
\endalign
$$
where
$$
J=(\sigma-\rho)A^2+2(-\sigma\alpha +\alpha -1)AB
+(\alpha^2\sigma +1-\alpha^2)B^2.
$$
Set $\sigma=(\alpha -1)/\alpha$. Then 
$$\align
J=&(\sigma-\rho)A^2-(\alpha-1)B^2\\
\ge&(\sigma-\rho)A^2-\frac{\alpha-1}{\alpha^2}A^2\\
=&\biggl(\frac{\alpha-1}{\alpha}-\frac{\alpha-1}{\alpha^2}-\rho\biggr )
A^2\\
=&\biggl(\biggl (\frac{\alpha-1}{\alpha}\biggr )^2-\rho\biggr )A^2\\
\ge&0\quad\forall 0<\rho<(\alpha-1)^2/\alpha^2.
\endalign
$$
Hence
$$
I((\alpha-1)^2/2\alpha^2)\ge(1-\sigma)(A-\alpha B)^2
=\frac{1}{\alpha}(A-\alpha B)^2.
$$

\noindent $\underline{\text{\bf Case 2}}$: $B\le 0$.

Let $\rho=(\alpha -1)^2/\alpha^2$. Then
$$
I(\rho)=\frac{1}{\alpha^2}[(A-\alpha B)+(\alpha -1)A]^2-\rho A^2
=\frac{1}{\alpha^2}(A-\alpha B)^2+2\frac{(\alpha -1)}{\alpha^2}A(A-\alpha B)
\ge\frac{1}{\alpha^2}(A-\alpha B)^2.
$$
By case 1 and case 2 the lemma follows.
\enddemo

\proclaim{\bf Theorem 1.2}
There exists a constant $R_0>0$ such that if $\1 M$ satisfies the 
interior rolling $R$-ball condition for some $0<R\le R_0$, then for any 
$\alpha>1+H$, $0<\delta<1$, there exists a constant $C_1>0$ depending on 
$k_0$, $H$, $\alpha$, $\delta$ and the space-time uniform bound of 
$|h_{ij}|$, $|\nabla^th_{ij}|$, $|q|$, $|\nabla^tq|$, 
$|\Delta^tq|$, such that
$$
\frac{|\nabla u|^2}{u^2}-\alpha\frac{u_t}{u}\le C_1
+\frac{\alpha^2(\alpha -1)^2(1+H)^2}{(1-\delta)^2(\alpha -(1+H))^2}\cdot
\frac{n}{2t}\quad\text{ in }M\times (0,T].\tag 1.5
$$
\endproclaim
\demo{Proof}
We will use a modification of the argument of \cite{CTY}, \cite{LY}, and 
\cite{W} to prove the theorem. Suppose $\1 M$ satisfies the interior rolling 
$R$-ball condition for some $0<R\le R_0$ where $R_0>0$ is some constant to 
be determined later. By \cite{C} there exists a  
$C^2$-function $\psi:[0,\infty)\to [0,H]$ such that $\psi (0)=0$,
$\psi (r)=H$ for all $r\ge 1$, which satisfy $0\le\psi'(r)\le 2H$ for all 
$r\ge 0$, $\psi'(0)=H$, and $\phi''(r)\ge -H$ for all $r\ge 0$. Let
$$
\phi (x)=1+\psi\biggl (\frac{r(x)}{R}\biggr )
$$
and $f=\log u$. By \cite{LY} and \cite{CTY}, $f$ satisfies
$$
\Delta f-f_t=q-|\nabla f|^2\quad\text{ in } M\times (0,T)\tag 1.6
$$
and
$$\left\{\aligned
(\Delta f)_t=&\Delta f_t-2h_{ij}f_{ij}-2h_{ik;i}f_k+\nabla (g^{ij}h_{ij})
\cdot\nabla f\\
(|\nabla f|^2)_t=&2\nabla f_t\cdot\nabla f-2h(\nabla f,\nabla f)
\endaligned\right.
$$
where $h(\nabla f,\nabla f)=h_{ij}f_if_j$. Let
$$
F(x,t)=t[\phi (x)(|\nabla f|^2+1)-\alpha f_t-\alpha q].
$$
Then in normal coordinates
$$\align
t^{-1}\Delta F=&\phi\cdot\Delta(|\nabla f|^2+1)+2\nabla\phi\cdot\nabla
|\nabla f|^2+\Delta\phi\cdot (|\nabla f|^2+1)-\alpha\Delta f_t
-\alpha\Delta q\\
=&2\phi (\sum_{i,j}f_{ij}^2+f_if_{ijj})+2\nabla\phi\cdot\nabla|\nabla f|^2
+\Delta\phi\cdot (|\nabla f|^2+1)-\alpha(\Delta f)_t-2\alpha h_{ij}f_{ij}\\
&\qquad -2\alpha h_{ik;i}f_k+\alpha\nabla (g^{ij}h_{ij})
\cdot\nabla f-\alpha\Delta q\\
=&2\biggl [\phi\sum_{i,j}f_{ij}^2-\alpha h_{ij}f_{ij}\biggr ]
+2\phi\nabla f\cdot\nabla(\Delta f)+2\phi R_{ij}f_if_j+2\nabla\phi\cdot\nabla
|\nabla f|^2\\
&\qquad +\Delta\phi\cdot (|\nabla f|^2+1)-\alpha(\Delta f)_t
-2\alpha h_{ik;i}f_k+\alpha\nabla (g^{ij}h_{ij})
\cdot\nabla f-\alpha\Delta q.\tag 1.7
\endalign
$$
By (1.6),
$$\align
&-\alpha(\Delta f)_t+2\phi\nabla f\cdot\nabla(\Delta f)\\
=&-\alpha(q+f_t-|\nabla f|^2)_t+2\phi\nabla f\cdot\nabla(q+f_t-|\nabla f|^2)\\
=&\alpha (-(\phi/\alpha)(|\nabla f|^2+1)+(F/(\alpha t))+|\nabla f|^2)_t
+2\phi\nabla f\cdot\nabla(q+f_t-|\nabla f|^2)\\
=&\frac{F_t}{t}-\frac{F}{t^2}+2\alpha\nabla f\cdot\nabla f_t
+2\phi\nabla f\cdot\nabla(q-|\nabla f|^2)+2(\phi-\alpha)h_{ij}f_if_j\\
=&\frac{F_t}{t}-\frac{F}{t^2}+2\alpha\nabla f\cdot\nabla\biggl (
\frac{\phi}{\alpha}(|\nabla f|^2+1)-\frac{F}{\alpha t}-q\biggr )
+2\phi\nabla f\cdot\nabla(q-|\nabla f|^2)+2(\phi-\alpha)h_{ij}f_if_j\\
=&\frac{F_t}{t}-\frac{F}{t^2}-\frac{2}{t}\nabla f\cdot\nabla F+2(\phi -\alpha)
\nabla f\cdot\nabla q-2\phi\nabla f\cdot\nabla |\nabla f|^2
+2\nabla f\cdot\nabla (\phi (|\nabla f|^2+1))\\
&\qquad +2(\phi-\alpha)h_{ij}f_if_j\\
=&\frac{F_t}{t}-\frac{F}{t^2}-\frac{2}{t}\nabla f\cdot\nabla F+2(\phi -\alpha)
\nabla f\cdot\nabla q+2(|\nabla f|^2+1)\nabla f\cdot\nabla\phi
+2(\phi-\alpha)h_{ij}f_if_j.\tag 1.8
\endalign
$$
Hence by (1.7) and (1.8),
$$\align
&t^{-1}(\Delta F-F_t+2\nabla f\cdot\nabla F)\\
=&2\biggl [\phi\sum_{i,j}f_{ij}^2-\alpha h_{ij}f_{ij}+2\phi_if_if_{ij}\biggr ]
+2(\phi-\alpha)\nabla f\cdot\nabla q+2\phi R_{ij}f_if_j
+2(\phi-\alpha)h_{ij}f_if_j-\frac{F}{t^2}\\
&\qquad +2(|\nabla f|^2+1)\nabla f\cdot\nabla\phi
+\Delta\phi\cdot (|\nabla f|^2+1)-2\alpha h_{ik;i}f_k+\alpha\nabla 
(g^{ij}h_{ij})\cdot\nabla f-\alpha\Delta q.\tag 1.9
\endalign
$$
By (1.9) and Young's inequality, for any $0<\delta<1$
there exist constants $C_1>0$, $C_2>0$, $C_3>0$, $C_4>0$, such that
$$\align
&t^{-1}(\Delta F-F_t+2\nabla f\cdot\nabla F)\\
\ge&2(\phi -\delta)\sum_{i,j}f_{ij}^2-C_1(|\nabla f|+|\nabla f|^2+|\nabla f|^3)
-C_2+\Delta\phi\cdot (|\nabla f|^2+1)-\frac{F}{t^2}\\
\ge&2(\phi -\delta)\sum_{i,j}f_{ij}^2-C_3|\nabla f|^3
-C_4+\Delta\phi\cdot (|\nabla f|^2+1)-\frac{F}{t^2}\\
\ge&\frac{2(\phi -\delta)}{n}(\Delta f)^2-C_3|\nabla f|^3
-C_4+\Delta\phi\cdot (|\nabla f|^2+1)-\frac{F}{t^2}.\tag 1.10
\endalign
$$
Let $R_1>0$ be the maximum number satisfying
$$\left\{\aligned
&\sqrt{k_0}\tan (R_1\sqrt{k_0})\le\frac{H}{2}+\frac{1}{2}\\
&\frac{H}{\sqrt{k_0}}\tan (R_1\sqrt{k_0})\le\frac{1}{2}.
\endaligned\right.
$$
and let $R_0\le R_1$. Then by the index comparison theorem (P.347 of \cite{Wa})
and an argument similar to that of \cite{C} there exists a constant
$c_1>0$ such that
$$
\nabla_i^0\nabla_j^0 r(x)\ge -c_1g_{ij}(x,0)\quad\forall x\in M, r(x)\le R_1.
\tag 1.11
$$
By (1.1) there exist constants $c_2>0$, $c_3>0$, such that
$$\left\{\aligned
&c_2g_{ij}(x,0)\le g_{ij}(x,t)\le c_3g_{ij}(x,0)\quad\forall 
x\in M, 0\le t\le T\\
&c_2g^{ij}(x,0)\le g^{ij}(x,t)\le c_3g^{ij}(x,0)\quad\forall 
x\in M, 0\le t\le T.\endaligned\right.\tag 1.12
$$
By (1.11), (1.12), and an argument similar to the proof of Lemma 1.3
of \cite{Hs1} there exists a constant $c_4>0$ such that
$$\align
\Delta^tr(x)\ge&-c_4\quad\forall x\in M, r(x)\le R_1, 0\le t\le T\\
\Rightarrow\qquad\Delta^t\phi=&\psi''\frac{|\nabla^tr(x)|^2}{R^2}+\psi'
\frac{\Delta^tr(x)}{R}
\ge-c_3\frac{H}{R^2}-c_4\frac{H}{R}\quad\forall x\in M, 0\le t\le T, 
r(x)\le R_1.\tag 1.13
\endalign
$$
Since $\Delta^t\phi=0$ for any $r(x)\ge R_1$ and $0\le t\le T$, by
(1.13) there exists a constant $c_5>0$ such that
$$
\Delta^t\phi\ge -c_5\quad\forall x\in M, 0\le t\le T.\tag 1.14
$$
By (1.10) and (1.14),
$$\align
t^{-1}(\Delta F-F_t+2\nabla f\cdot\nabla F)
\ge&\frac{2(\phi -\delta)}{n}(\Delta f)^2-C_3|\nabla f|^3
-C_4-C_5|\nabla f|^2-\frac{F}{t^2}\\
\ge&\frac{2(\phi -\delta)}{n}(\Delta f)^2-C_6|\nabla f|^3
-C_7-\frac{F}{t^2}.\tag 1.15
\endalign
$$
for some constants $C_5\ge 0$, $C_6>0$, $C_7>0$.

Since $M\times [0,T]$ is compact, there exists $(x_0,t_0)\in M\times [0,T]$
such that $F(x_0,t_0)=\max_{M\times [0,T]}F$. If $t_0=0$, $F\le 0$
on $M\times [0,T]$. Then (1.5) holds and we are done. Hence we may assume 
without loss of generality that $t_0>0$. Suppose first $x_0\in\1 M$. Since 
$f_{\nu}=u_{\nu}/u=0$ on $\1 M\times (0,T)$, we have
$$\align
0\le&\frac{1}{t_0}\frac{\1 F}{\1\nu}(x_0,t_0)=(|\nabla f|^2+1)\phi_{\nu}
+2\phi f_if_{i\nu}-\alpha f_{\nu t}-\alpha q_{\nu}\\
=&(|\nabla f|^2+1)\phi_{\nu}-2II(\nabla f,\nabla f)-\alpha q_{\nu}\\
\le&-\frac{H}{R}(|\nabla f|^2+1)+2H|\nabla f|^2+\alpha a_0\\
\le&H|\nabla f|^2(2-(1/R_0))+\alpha a_0-(H/R_0)\tag 1.16
\endalign
$$
where $a_0=\max_{M\times [0,T]}|\nabla^tq|$. Let 
$$
R_0=\min (1/2,R_1,H/(1+\alpha a_0)).
$$
Then the right hand side of (1.16) is strictly less than $0$. Hence 
contradiction arises. Thus $x_0\in M\setminus\1 M$. Then $\nabla F(x_0,t_0)=0$,
$\Delta F(x_0,t_0)\le 0$ and $F_t(x_0,t_0)\ge 0$. Hence at $(x_0,t_0)$,
$$
\Delta F-F_t+2\nabla f\cdot\nabla F\le 0.\tag 1.17
$$ 
By (1.6), (1.15), and (1.17), at $(x_0,t_0)$,
$$
\frac{2(\phi -\delta)}{n}(|\nabla f|^2-f_t-q)^2-C_6|\nabla f|^3
-C_7-\frac{F}{t_0^2}\le 0.\tag 1.18
$$
By an argument similar to that of P.382 of \cite{W},
$$
(|\nabla f|^2-f_t-q)^2\ge (1-\delta)(|\nabla f|^2+1-f_t-q)^2-\frac{2}{\delta}.
\tag 1.19
$$
If $F(x_0,t_0)<0$, then $\max_{M\times [0,T]}F<0$. Then (1.5) holds and we 
are done. Hence we may assume without loss of generality that $F(x_0,t_0)
\ge 0$. Then
$$
f_t+q\le\frac{\phi}{\alpha}(|\nabla f|^2+1)\tag 1.20
$$
at $(x_0,t_0)$. Let
$$
\delta_1=\frac{\alpha -(1+H)}{(\alpha -1)(1+H)}.
$$
Then 
$$
\frac{1+\delta_1(\alpha -1)}{\alpha}=\frac{1}{1+H}.\tag 1.21
$$
Hence by (1.20) and (1.21),
$$\align
&(|\nabla f|^2+1-f_t-q)^2-\delta_1^2(\phi\cdot(|\nabla f|^2+1)-f_t-q)^2\\
=&(|\nabla f|^2+1-f_t-q+\delta_1(\phi\cdot(|\nabla f|^2+1)-f_t-q))\\
&\qquad\cdot(|\nabla f|^2+1-f_t-q-\delta_1(\phi\cdot(|\nabla f|^2+1)-f_t-q))\\
=&((1+\delta_1\phi)(|\nabla f|^2+1)-(1+\delta_1)(f_t+q))
((1-\delta_1\phi)(|\nabla f|^2+1)-(1-\delta_1)(f_t+q))\\
\ge&((1+\delta_1\phi)-(1+\delta_1)(\phi/\alpha))
((1-\delta_1\phi)-(1-\delta_1)(\phi/\alpha))
(|\nabla f|^2+1)^2\\
=&\biggl (1+\delta_1\biggl (1-\frac{1}{\alpha}\biggr )\phi 
-\frac{\phi}{\alpha}\biggr )\biggl (
1-\frac{1+\delta_1(\alpha -1)}{\alpha}\phi\biggr )\\
\ge&\biggl (1-\frac{\phi}{\alpha}\biggr )\biggl (1-\frac{\phi}{1+H}\biggr )\\
\ge&\biggl (1-\frac{\phi}{1+H}\biggr )^2\\
\ge&0.\tag 1.22
\endalign
$$
By (1.18), (1.19) and (1.22), for any $0<\rho<1$ there exists a constant
$C_8=C_8(\rho)>0$ such that
$$
0\ge\frac{2(1-\delta)^2(\alpha -(1+H))^2}{n(\alpha -1)^2(1+H)^2}
\{(\phi\cdot(|\nabla f|^2+1)-f_t-q)^2-\rho (\phi\cdot(|\nabla f|^2+1))^2\}
-C_8-\frac{F}{t_0^2}.\tag 1.23
$$
Let $A=\phi\cdot(|\nabla f|^2+1)$ and $B=f_t+q$. By Lemma 1.1 there exists a 
constant $0<\rho<1$ such that (1.4) holds. By (1.4) and (1.23),
$$\align
&\frac{2(1-\delta)^2(\alpha -(1+H))^2}{n\alpha^2(\alpha -1)^2(1+H)^2}
(\phi\cdot(|\nabla f|^2+1)-\alpha (f_t+q))^2-C_8-\frac{F}{t_0^2}\le 0\\
\Rightarrow\quad
&\frac{2(1-\delta)^2(\alpha -(1+H))^2}{n\alpha^2(\alpha -1)^2(1+H)^2}F^2
-F-C_8t_0^2\le 0\\
\Rightarrow\quad&\biggl (
F-\frac{n\alpha^2(\alpha -1)^2(1+H)^2}{4(1-\delta)^2(\alpha -(1+H))^2}
\biggr )^2\le\biggl (\frac{n\alpha^2(\alpha -1)^2(1+H)^2}{4(1-\delta)^2
(\alpha -(1+H))^2}\biggr )^2+C_9t_0^2\\
\Rightarrow\quad&F(x_0,t_0)\le
\frac{n\alpha^2(\alpha -1)^2(1+H)^2}{2(1-\delta)^2(\alpha -(1+H))^2}
+C_{10}t_0
\le\frac{n\alpha^2(\alpha -1)^2(1+H)^2}{2(1-\delta)^2(\alpha -(1+H))^2}
+C_{10}T\\
\Rightarrow\quad&F(x,T)\le
\frac{n\alpha^2(\alpha -1)^2(1+H)^2}{2(1-\delta)^2(\alpha -(1+H))^2}+C_{10}T
\quad\forall x\in M.\tag 1.24
\endalign
$$
By replacing $T$ by $t$ in (1.24) for any $t\in (0,T]$ the theorem follows.
\enddemo

\proclaim{\bf Corollary 1.3}
Suppose $\partial M$ is convex with respect to $g(t)$ for all $0\le t\le T$.
Then for any $0<\3<1$ there exists a constant 
$C_1>0$ depending on $k_0$, $\3$, and the space-time uniform bound of 
$|h_{ij}|$, $|\nabla^th_{ij}|$, $|q|$, $|\nabla^tq|$, 
$|\Delta^tq|$, such that
$$
\frac{|\nabla u|^2}{u^2}-(1+\3+\3^2)\frac{u_t}{u}\le C_1
+\frac{(1+\3+\3^2)^2(1+\3)^2(1+\3^2)^2}{1-\3^2}\cdot
\frac{n}{2t}\quad\text{ in }M\times (0,T].\tag 1.25
$$
\endproclaim
\demo{Proof}
(1.25) follows from (1.5) by setting $\alpha=1+\3+\3^2$ and $\delta=H=\3^2$.
\enddemo

By an argument similar to the proof of Corollary 4.2 of \cite{CTY} but with
Theorem 1.2 replacing Lemma 4.1 of \cite{CTY} in the proof there we get
the following corollary.

\proclaim{\bf Corollary 1.4}
Let $R_0>0$ be given by Theorem 1.2. Suppose $\1 M$ satisfies the 
interior rolling $R$-ball condition for some $0<R\le R_0$.
Then for any $\alpha >1+H$, $0<\delta<1$, there exist constants 
$C_2>0$, $C_3>0$, depending on $k_0$, $H$, $\alpha$, $\delta$ and the 
space-time uniform bound of $|h_{ij}|$, $|\nabla^th_{ij}|$, $|q|$, 
$|\nabla^tq|$, $|\Delta^tq|$, such that
$$
u(x_1,t_1)\le u(x_2,t_2)\biggl (\frac{t_2}{t_1}\biggr )^{a_1}\text{exp }
\biggl (C_2\alpha\frac{r(x_1,x_2)}{t_2-t_1}+C_3(t_2-t_1)\biggr )
\tag 1.26
$$
for any $x_1,x_2\in M$, $0<t_1<t_2\le T$, where
$$
a_1=\frac{n\alpha^2(\alpha -1)^2(1+H)^2}{2(1-\delta)^2(\alpha -(1+H))^2}.
$$
\endproclaim

By an argument similar to the proof of Theorem 1.2 we have the following
theorem.

\proclaim{\bf Theorem 1.5}
Suppose $(M,g(t))$, $0\le t\le T$, is a closed manifold with metric $g(t)$
satisfying (1.1) and (1.2) for some smooth symmetric tensors $h_{ij}(x,t)$ 
on $M\times [0,T]$ and constant $k_0>0$. Let $u$ be the solution of 
$$
u_t=\Delta^tu-qu\quad\text{ in }M\times [0,T]
$$
for some smooth function $q(x,t)$ on $M\times [0,T]$. Then for any $\alpha>1$ 
there exist constants $C_1>0$, $C_2>0$, $C_3>0$, depending on $k_0$,
$\alpha$, and the space-time uniform bound of $|h_{ij}|$, 
$|\nabla^th_{ij}|$, $|q|$, $|\nabla^tq|$, $|\Delta^tq|$, such that 
(1.5) and (1.26) holds with $\delta=H=0$.
\endproclaim

$$
\text{Section 2}
$$

In this section we will give another proof of Perelman's Li-Yau-Hamilton type 
inequality for the fundamental solution of the conjugate heat equation on 
closed manifolds without using the properties of the reduced distance. 
Let $(M,g(t))$, $0\le t\le T$, be a closed manifold with metric $g(t)$ 
evolving by the Ricci flow (0.1). 
Let $\Cal{Z}(x,t;y,s)$, $0\le s<t\le T$, be the heat kernel of $M$. Then
$\forall y\in M$, $0\le s<T$, $\Cal{Z}(\cdot,\cdot;y,s)$ satisfies
$$\left\{\aligned
&\Cal{Z}_t=\Delta^t_x\Cal{Z}\quad\text{ in }M\times (s,T)\\
&\lim_{t\to s}\int_M\Cal{Z}(z,t,y,s)\eta (z)\,dV_t(z)=\eta (y)\quad\forall
\eta\in C^{\infty}(M)
\endaligned\right.
$$
and for any $x\in M$, $0<t\le T$, $\Cal{Z}(x,t;\cdot,\cdot)$ satisfies
$$\left\{\aligned
&-\Cal{Z}_s-\Delta_y\Cal{Z}+R\Cal{Z}=0\quad\text{ in }M\times (0,t)\\
&\lim_{s\to t}\int_M\Cal{Z}(x,t,y,s)\eta (y)\,dV_t(y)=\eta (x)\quad\forall
\eta\in C^{\infty}(M).
\endaligned\right.
$$
By an argument similar to the proof of Corollary 5.2 of \cite{CTY} we
have the following result.

\proclaim{\bf Lemma 2.1}
There exist constants $C>0$ and $D>0$ such that
$$\left\{\aligned
&\Cal{Z}(x,t;y,s)\le\frac{C}{V_x(\sqrt{t-s})}
\text{exp }\biggl (-\frac{r^2(x,y)}{D(t-s)}\biggr )
\quad\forall x,y\in M,0\le s<t\le T\\
&\Cal{Z}(x,t;y,s)\le\frac{C}{V_y(\sqrt{t-s})}
\text{exp }\biggl (-\frac{r^2(x,y)}{D(t-s)}\biggr )
\quad\forall x,y\in M,0\le s<t\le T.\endaligned\right.
$$ 
\endproclaim

Let $p\in M$ and $u(x,t)=\Cal{Z}(p,T,x,t)$. Then $u$ satisfies (0.3) and 
(0.4) in $M\times (0,T)$. As in \cite{P} we let $f$, $v$, be 
given by (0.2) and (0.5) with $\tau=T-t$. Let $0\le h_0\in C^{\infty}(M)$, 
$0<t_0<T$, and let $h\ge 0$ be the solution of the heat equation
$$\left\{\aligned
h_t=&\Delta h\qquad\text{ in }M\times (t_0,T]\\
h(x,t_0)=&h_0(x)\quad\text{ in }M.
\endaligned\right.
$$
We next recall a result of Perelman \cite{P}.

\proclaim{\bf Lemma 2.2}(\cite{P})
$$
\int_Mvh\,dV_{t_1}\le\int_Mvh\,dV_{t_2}\quad\forall t_0\le t_1\le t_2<T.
$$
\endproclaim

\proclaim{\bf Lemma 2.3}
$$
\limsup_{t\to T}\int_M\tau hu(2\Delta f-|\nabla f|^2+R)\,dV_t
\le\frac{n}{2}h(p,T)\tag 2.1
$$
where $\tau=T-t$.
\endproclaim
\demo{Proof}
We will use a modification of the technique of \cite{CTY} to prove the lemma.
By direct computation,
$$\align
&\limsup_{t\to T}\int_M\tau hu(2\Delta f-|\nabla f|^2+R)\,dV_t\\
=&\limsup_{t\to T}\int_M\tau h\biggl (-2\Delta u+\frac{|\nabla u|^2}{u}
+Ru\biggr )\,dV_t\\
=&-2\lim_{t\to T}\tau\int_Mu\Delta h\,dV_t
+\limsup_{t\to T}\tau\int_Mh\frac{|\nabla u|^2}{u}\,dV_t
+\lim_{t\to T}\tau\int_MRhu\,dV_t.
\endalign
$$
Note that
$$
\biggl |\tau\int_MRhu\,dV_t\biggr |\le\tau\|R\|_{\infty}\|h\|_{\infty}\to
0\quad\text{ as }t\to T.
$$
Since by the Schauder estimates \cite{LSU},
$$
\sup_{\frac{T+t_0}{2}\le s\le T}\|\Delta h(\cdot,s)\|_{L^{\infty}(M)}<\infty,
$$
$$
\biggl |\tau\int_Mu\Delta h\,dV_t\biggr |\le\tau
\sup_{\frac{T+t_0}{2}\le s\le T}\|\Delta h(\cdot,s)\|_{L^{\infty}(M)}
\to 0\quad\text{ as }t\to T.
$$
Hence
$$
\limsup_{t\to T}\int_M\tau hu(2\Delta f-|\nabla f|^2+R)\,dV_t
=\limsup_{t\to T}\int_M\tau h\frac{|\nabla u|^2}{u}\,dV_t.\tag 2.2
$$
By Theorem 1.5 for any $\alpha>1$ there exists a constant $C_1>0$ such that
$$
\frac{|\nabla^tu|^2}{u^2}-\alpha\frac{u_{\tau}}{u}\le C_1
+\frac{n\alpha^2}{2\tau}\quad\text{ in }M\times (0,T]
$$
where $\tau=T-t$. Then
$$\align
\tau\int_Mh\frac{|\nabla u|^2}{u}\,dV_t
\le&\tau\int_Mh\biggl (\alpha u_{\tau}+C_1u+\frac{n\alpha^2}{2\tau}u\biggr )
\,dV_t\\
=&\tau\int_M[\alpha h(\Delta u-Ru)+C_1hu]\,dV_t
+\frac{n\alpha^2}{2}\int_Mhu\,dV_t\\
=&\tau\int_M[\alpha (u\Delta h-Ru)+C_1hu]\,dV_t
+\frac{n\alpha^2}{2}\int_Mhu\,dV_t.\tag 2.3
\endalign
$$
Since
$$
\biggl |\tau\int_M[\alpha(u\Delta h-Ru)+C_1hu]\,dV_t\biggr |\le C\tau\to 0
\quad\text{ as }\tau\to 0
$$
and
$$
\lim_{t\to T}\int_Mhu\,dV_t=h(p,T),\tag 2.4
$$
letting $t\to T$ in (2.3) we get
$$\align
&\limsup_{t\to T}\tau\int_Mh\frac{|\nabla u|^2}{u}\,dV_t\le
\frac{n\alpha^2}{2}h(p,T)\quad\forall\alpha>1\\
\Rightarrow\quad&\limsup_{t\to T}\tau\int_Mh\frac{|\nabla u|^2}{u}\,dV_t\le
\frac{n}{2}h(p,T)\quad\text{ as }\alpha\to 1.\tag 2.5
\endalign
$$
By (2.2) and (2.5) we get (2.1) and the lemma follows. 
\enddemo

By the same argument as the proof of Lemma 7.6 of \cite{CTY} but with 
Lemma 2.1 replacing Corollary 5.2 of \cite{CTY} in the proof there we get

\proclaim{\bf Lemma 2.4}
$$
\limsup_{t\to T}\int_Mfhu\,dV_t\le \frac{n}{2}h(p,T).
$$
\endproclaim

By (2.4), Lemma 2.3, and Lemma 2.4 we get

\proclaim{\bf Lemma 2.5}
$$
\limsup_{t\to T}\int_Mvh\,dV_t\le 0.
$$
\endproclaim

\proclaim{\bf Theorem 2.6}(Corollary 9.3 of \cite{P})
(0.6) holds in $M\times (0,T)$.
\endproclaim
\demo{Proof}
By Lemma 2.2 and Lemma 2.5,
$$
\int_Mv(x,t_0)h_0(x)\,dV_{t_0}(x)\le 0\quad\forall 0\le h_0\in 
C^{\infty}(M), 0<t_0<T 
$$
and the theorem follows.
\enddemo

$$
\text{ Section 3}
$$

In this section we will generalize a result of Q.S.~Zhang \cite{Z} to local
gradient estimates for the solution of the generalized conjugate heat 
equation on compact manifolds.

\proclaim{\bf Theorem 3.1}
Let $(M,g(t))$, $0\le t\le T$, be a $n$-dimensional compact manifold, $n\ge 2$,
with metric $g(t)$ satisfying (1.1) and (1.2) for some smooth family of 
symmetric tensors $h_{ij}(x,t)$ on $M$ and constant $k_0>0$. Let $u$ be 
a positive solution of 
$$
u_t=\Delta^tu-qu\quad\text{ in }M\times [0,T] 
$$
where $q(x,t)$ is a smooth function on $M\times [0,T]$. Let $x_0\in M$ and
$t_0\in (0,T]$. Then there exists a constant $C_1>0$  depending on $k_0$ and 
the space-time uniform bound of $|h_{ij}|$, $|q|$ and $|\nabla^tq|$ in
$Q_{R,T_1}(x_0,t_0)$ such that
$$
\frac{|\nabla u|}{u}\le C_1\biggl (\frac{1}{R}+\frac{1}{\sqrt{T_1}}+1
\biggr )\biggl (1+\log\biggl(\frac{A}{u}\biggr )\biggr )
$$
holds in $Q_{R/2,T_1/2}(x_0,t_0)$ for any $Q_{R,T_1}=Q_{R,T_1}(x_0,t_0)
\subset M\setminus\1 M$, $0<R\le 1$, $0<T_1\le t_0$, whenever $u\le A$ in 
$Q_{R,T_1}(x_0,t_0)$ for some constant $A>0$. 
\endproclaim
\demo{Proof}
We will use a modification of the proof of Theorem 3.1 of \cite{Z} 
and Theorem 1.1 of \cite{SZ} to prove the theorem. Suppose $Q_{R,T_1}\subset 
M\setminus\1 M$, $0<R\le 1$, and $0<T_1\le t_0$.
Since (3.1) is invariant by rescaling $u$ to $u/A$, we may 
assume without loss of generality that $0<u\le 1$ in $Q_{R,T_1}(x_0,t_0)$.
As in \cite{Z} let $f=\log u$ and 
$$
w=\frac{|\nabla f|^2}{(1-f)^2}.
$$
Since $f$ satisfies (1.6), by a direct computation we have in normal 
coordinates,
$$\align
w_t=&2\frac{\nabla f\cdot\nabla f_t}{(1-f)^2}
+2\frac{|\nabla f|^2f_t}{{(1-f)^3}}-2\frac{h_{ij}f_if_j}{(1-f)^2}\\
=&2\frac{\nabla f\cdot\nabla (\Delta f+|\nabla f|^2-q)}{(1-f)^2}
+2\frac{|\nabla f|^2(\Delta f+|\nabla f|^2-q)}{(1-f)^3}
-2\frac{h_{ij}f_if_j}{(1-f)^2}\tag 3.1\\
\nabla_jw=&2\frac{f_if_{ij}}{(1-f)^2}+2\frac{|\nabla f|^2f_j}{(1-f)^3}
\tag 3.2\\
\Delta w=&\frac{2f_{ij}^2}{(1-f)^2}+2\frac{f_if_{ijj}}{(1-f)^2}
+8\frac{f_if_jf_{ij}}{(1-f)^3}+2\frac{|\nabla f|^2\Delta f}{(1-f)^3}
+6\frac{|\nabla f|^4}{(1-f)^4}.\tag 3.3
\endalign
$$
By (3.1), (3.2), and (3.3),
$$\align
&\Delta w-w_t\\
=&\frac{2f_{ij}^2}{(1-f)^2}+8\frac{f_if_jf_{ij}}{(1-f)^3}
+2\frac{f_i(f_{ijj}-f_{jji})}{(1-f)^2}
-2\frac{\nabla f\cdot\nabla (|\nabla f|^2-q)}{(1-f)^2}
+2\frac{h_{ij}f_if_j}{(1-f)^2}\\
&\qquad -2\frac{|\nabla f|^2(|\nabla f|^2-q)}{(1-f)^3}
+6\frac{|\nabla f|^4}{(1-f)^4}\\
=&\frac{2f_{ij}^2}{(1-f)^2}+8\frac{f_if_jf_{ij}}{(1-f)^3}
+2\frac{R_{ij}f_if_j}{(1-f)^2}
-\frac{2\nabla f\cdot\nabla (|\nabla f|^2-q)}{(1-f)^2}
+2\frac{h_{ij}f_if_j}{(1-f)^2}\\
&\qquad -2\frac{|\nabla f|^2(|\nabla f|^2-q)}{(1-f)^3}
+6\frac{|\nabla f|^4}{(1-f)^4}\\
=&\frac{2f_{ij}^2}{(1-f)^2}+8\frac{f_if_jf_{ij}}{(1-f)^3}
+6\frac{|\nabla f|^4}{(1-f)^4}-\frac{4f_if_jf_{ij}}{(1-f)^2}
-2\frac{|\nabla f|^4}{(1-f)^3}
+2\frac{h_{ij}f_if_j}{(1-f)^2}\\
&\qquad +2\frac{R_{ij}f_if_j}{(1-f)^2}
+2\frac{|\nabla f|^2q}{(1-f)^3}
+\frac{2\nabla f\cdot\nabla q}{(1-f)^2}\\
=&\frac{2}{(1-f)^2}\biggl (f_{ij}+\frac{f_if_j}{1-f}\biggr )^2
+4\frac{f_if_jf_{ij}}{(1-f)^3}
+4\frac{|\nabla f|^4}{(1-f)^4}-\frac{4f_if_jf_{ij}}{(1-f)^2}
-2\frac{|\nabla f|^4}{(1-f)^3}\\
&\qquad +2\frac{h_{ij}f_if_j}{(1-f)^2}
+2\frac{R_{ij}f_if_j}{(1-f)^2}
+2\frac{|\nabla f|^2q}{(1-f)^3}
+\frac{2\nabla f\cdot\nabla q}{(1-f)^2}.
\endalign
$$
Hence
$$\align
&\Delta w-w_t\\
\ge&4\frac{f_if_jf_{ij}}{(1-f)^3}
+4\frac{|\nabla f|^4}{(1-f)^4}-\frac{4f_if_jf_{ij}}{(1-f)^2}
-2\frac{|\nabla f|^4}{(1-f)^3}+2\frac{h_{ij}f_if_j}{(1-f)^2}
+2\frac{R_{ij}f_if_j}{(1-f)^2}\\
&\qquad +2\frac{|\nabla f|^2q}{(1-f)^3}
+\frac{2\nabla f\cdot\nabla q}{(1-f)^2}\\
=&\frac{2}{(1-f)}\biggl (\nabla f\cdot\nabla w
-2\frac{|\nabla f|^4}{(1-f)^3}\biggr )
+4\frac{|\nabla f|^4}{(1-f)^4}+2\biggl (2\frac{|\nabla f|^4}{(1-f)^3}
-\nabla f\cdot\nabla w\biggr )\\
&\qquad-2\frac{|\nabla f|^4}{(1-f)^3}+2\frac{h_{ij}f_if_j}{(1-f)^2}
+2\frac{R_{ij}f_if_j}{(1-f)^2}+2\frac{|\nabla f|^2q}{(1-f)^3}
+\frac{2\nabla f\cdot\nabla q}{(1-f)^2}\\
=&2\frac{|\nabla f|^4}{(1-f)^3}+\frac{2f}{1-f}\nabla f\cdot\nabla w
+2\frac{h_{ij}f_if_j}{(1-f)^2}+2\frac{R_{ij}f_if_j}{(1-f)^2}
+2\frac{|\nabla f|^2q}{(1-f)^3}
+\frac{2\nabla f\cdot\nabla q}{(1-f)^2}.\tag 3.4
\endalign
$$
Since $f\le 0$, by (3.4) there exist constants $C_1>0$ and $C_2>0$ 
depending on $k_0$ and the space-time uniform bound of $|h_{ij}|$, $|q|$
and $|\nabla^tq|$ in $Q_{R,T_1}(x_0,t_0)$ such that
$$
\Delta w-w_t\ge\frac{2f}{1-f}\nabla f\cdot\nabla w+2(1-f)w^2-C_1w-C_2
\quad\text{ in }Q_{R_1,T_1}.\tag 3.5
$$
We choose a smooth nonnegative function $\phi:\Bbb{R}\to\Bbb{R}$,
$0\le\phi\le 1$, such that $\phi (r)=1$ for all $r\le 1/2$, $\phi (r)=0$
for all $r\ge 1$, and $\phi'(r)\le 0$ for all $r\in\Bbb{R}$. For any $x\in M$, 
$0\le t\le T$, let $\phi_R(x)=(\phi(r(x,x_0)/R))^5$, $\eta_{T_1}(t)
=(\phi ((t_0-t)/T_1))^5$, and
$$
\psi (x,t)=\phi_R(x)\eta_{T_1}(t).
$$
When there is no ambiguity, we will write $r$ for $r(x,x_0)$. Similar to the 
proof of Theorem 1.2 by (1.1) there exist constants $c_2>0$ and $c_3>0$ such
that (1.12) holds in $Q_{R,T_1}(x_0,t_0)$. Then by (1.12),
$$
|\nabla^tr|^2\le c_3|\nabla^0r|^2\le c_3\quad\text{ in }Q_{R,T_1}(x_0,t_0).
$$
Hence
$$\align
\nabla^t\phi_R(x)=5\phi^4\phi'\cdot\frac{\nabla^tr}{R}\quad
\Rightarrow&\quad\frac{|\nabla^t\phi_R|}{\phi_R^{3/4}}=
5\frac{\phi^{\frac{1}{4}}|\phi'||\nabla^tr|}{R}\le\frac{C_3}{R}\\
\Rightarrow&\quad\frac{|\nabla^t\phi_R|}{\phi_R^{1/2}}\le\frac{C_3}{R}
\tag 3.6
\endalign
$$
for some constant $C_3>0$ and
$$
\Delta^t\phi_R(x)=(20\phi^3\phi'{}^2+5\phi^4\phi'')\frac{|\nabla^t r|^2}{R^2}
+5\phi^4\phi'\cdot\frac{\Delta^t r}{R}
\ge-\frac{C_3'}{R^2}
+5\phi^4\phi'\cdot\frac{\Delta^t r}{R}\tag 3.7
$$
for some constant $C_3'>0$. Similarly
$$
\frac{|\1_t\eta_{T_1}|}{\eta_{T_1}^{1/2}}\le\frac{C_4}{T_1}\tag 3.8
$$
for some constant $C_4>0$. By (1.2) and the Hessian comparison theorem 
\cite{SY},
$$
\nabla_i^0\nabla_j^0r\le\frac{n-1}{r}(1+\sqrt{k_0}r)g_{ij}(x,0).\tag 3.9
$$
By (1.1), (3.9), and an argument similar to the proof of Lemma 1.3 of 
\cite{Hs1} there exists a constant $C_4'>0$ such that
$$
\Delta^tr\le C_4'(1+\frac{1}{r})\quad\forall 0\le t\le T.\tag 3.10
$$
By (3.7) and (3.10),
$$
\frac{\Delta^t\phi_R}{\phi_R^{1/2}}
\ge-\frac{C_3'}{R^2}-5C_4'\phi^{3/2}|\phi'|(1+(R/2)^{-1})\ge-\frac{C_5}{R^2}
\tag 3.11
$$
for some constant $C_5>0$. By (3.6), (3.8) and (3.11), there exists a constant
$C_6>0$ such that
$$
\frac{|\nabla\psi|}{\psi^{1/2}}\le\frac{C_6}{R},
\frac{|\nabla\psi|}{\psi^{3/4}}\le\frac{C_6}{R},
\frac{|\1_t\psi|}{\psi^{1/2}}\le\frac{C_6}{T_1},
\frac{\Delta\psi}{\psi^{1/2}}\ge-\frac{C_6}{R^2}.\tag 3.12
$$
By (3.5),
$$\align
&\Delta (\psi w)-(\psi w)_t\\
=&\psi (\Delta w-w_t)+2\nabla\psi\cdot\nabla w+w\Delta\psi-\psi_tw\\
\ge&\frac{2f\psi}{1-f}\nabla f\cdot\nabla w+2(1-f)\psi w^2-C_1\psi w
-C_2\psi+2\nabla\psi\cdot\nabla w+w\Delta\psi-\psi_tw\\
\ge&\frac{2f}{1-f}\nabla f\cdot(\nabla (\psi w)-w\nabla\psi)
+2(1-f)\psi w^2-C_1\psi w-C_2\psi+2\frac{\nabla\psi}{\psi}\cdot
(\nabla (\psi w)-w\nabla\psi)\\
&\qquad +w\Delta\psi-\psi_tw\\
\ge&\frac{2f}{1-f}\nabla f\cdot\nabla (\psi w)
+2\frac{\nabla\psi}{\psi}\cdot\nabla (\psi w)-\frac{2fw}{1-f}\nabla f
\cdot\nabla\psi+2(1-f)\psi w^2-C_1\psi^{\frac{1}{2}}w-C_2\\
&\qquad -2w\frac{|\nabla\psi|^2}{\psi}+w\Delta\psi-\psi_tw.\tag 3.13
\endalign
$$
By (3.12),
$$\left\{\aligned
&w\frac{|\nabla\psi|^2}{\psi}=(\psi^{1/2}w)\biggl (
\frac{|\nabla\psi|}{\psi^{3/4}}\biggr )^2\le\frac{C_6^2}{R^2}\psi^{1/2}w
\le\frac{1}{16}\psi w^2+\frac{4C_6^4}{R^4}\\
&w\Delta\psi=(\psi^{1/2}w)\frac{\Delta\psi}{\psi^{1/2}}\ge-C_6
\frac{\psi^{1/2}w}{R^2}\ge -\frac{1}{8}\psi w^2-\frac{2C_6^2}{R^4}\\
&|w\1_t\psi|=(\psi^{1/2}w)\frac{|\1_t\psi|}{\psi^{1/2}}\le C_6
\frac{\psi^{1/2}w}{T_1}\le\frac{1}{8}\psi w^2+\frac{2C_6^2}{T_1^2}.
\endaligned\right.\tag 3.14
$$
Similarly,
$$
C_1\psi^{1/2}w\le\frac{1}{8}\psi w^2+2C_1^2.\tag 3.15
$$
By (3.12) and an argument similar to the proof of Theorem 3.1 of \cite{Z},
$$
\biggl |\frac{2fw}{1-f}\nabla f\cdot\nabla\psi\biggr |
\le (1-f)\psi w^2+C_7\frac{f^4}{R^4(1-f)^3}\tag 3.16
$$
for some constant $C_7>0$. Since
$$
\frac{|f|}{(1-f)}\le 1,
$$
by (3.13), (3.14), (3.15) and (3.16),
$$\align
&\Delta (\psi w)-(\psi w)_t-\frac{2f}{1-f}\nabla f\cdot\nabla (\psi w)
-2\frac{\nabla\psi}{\psi}\cdot\nabla (\psi w)\\
\ge&(1-f)\psi w^2-\frac{1}{2}\psi w^2-C_8\bigg (
\frac{1}{R^4}+\frac{1}{T_1^2}+1\biggr )-C_8\frac{f^4}{R^4(1-f)^3}\\
\ge&(1-f)\psi w^2-\frac{1}{2}\psi w^2-C_8\bigg (
\frac{1}{R^4}+\frac{1}{T_1^2}+1\biggr )-C_8\frac{(1-f)}{R^4}\tag 3.17
\endalign
$$
for some constant $C_8>0$. Suppose the function $\psi w$ attains its maximum
on the set $Q_{R,T_1}$ at the point $(x_1,t_1)\in\2{Q}_{R,T_1}$. 
Similar to \cite{LY} and \cite{Z} we may assume without loss of generality
that $x_1$ is not a cut point of $x_0$ with respect to the metric $g(0)$. 
Then at $(x_1,t_1)$, $\nabla (\psi w)=0$, $\1_t(\psi w)\ge 0$, 
$\Delta (\psi w)\le 0$. Hence the left hand side of (3.17) is $\le 0$ at 
$(x_1,t_1)$. Thus
$$\align
&(1-f)\psi w^2\le\frac{1}{2}\psi w^2+C_8\bigg (\frac{1}{R^4}
+\frac{1}{T_1^2}+1\biggr )+C_8\frac{(1-f)}{R^4}\\
\Rightarrow\quad&\frac{1}{2}(1-f)\psi w^2\le C_8\bigg (\frac{1}{R^4}
+\frac{1}{T_1^2}+1\biggr )+C_8\frac{(1-f)}{R^4}\\
\Rightarrow\quad&\sup_{\2{Q}_{R,T_1}}\psi w^2\le\psi(x_1,t_1)w^2(x_1,t_1)
\le 2C_8\bigg (\frac{2}{R^4}+\frac{1}{T_1^2}+1\biggr ).
\endalign
$$
Hence
$$\align
&w^2(x,t)\le 4C_8\bigg (\frac{1}{R^4}+\frac{1}{T_1^2}+1\biggr )
\qquad\qquad\qquad\qquad\quad\text{in }Q_{R/2,T_1/2}\\
\Rightarrow\quad&\frac{|\nabla f(x,t)|^2}{(1-f(x,t))^2}
=w(x,t)\le 2\sqrt{C_8}\bigg (\frac{1}{R^2}+\frac{1}{T_1}+
1\biggr )\quad\text{ in }Q_{R/2,T_1/2}\\
\Rightarrow\quad&\frac{|\nabla f(x,t)|}{(1-f(x,t))}
\le 2C_8^{1/4}\bigg (\frac{1}{R}+\frac{1}{\sqrt{T_1}}+1\biggr )
\qquad\qquad\qquad\text{in }Q_{R/2,T_1/2}
\endalign
$$
and the theorem follows.
\enddemo

$$
\text{ Section 4}
$$

In this section we will prove various gradient estimates for the Dirichlet 
fundamental solution of the conjugate heat equation. 

Let $(M,g(t))$, $0\le t\le T$, be a complete 
noncompact $n$-dimensional manifold, $n\ge 2$, with metric $g(t)$ evolving 
by the Ricci flow (0.1) which satisfies
$$
|\nabla^iRm|\le k_0\quad\text{ on }M\times [0,T]\quad\forall i=0,1,2\tag 4.1
$$
for some constant $k_0\ge 1$. Similar to section 2 we let $\Cal{Z}(x,t;y,s)$, 
$0\le s<t\le T$, be the fundamental solution of the heat equation in 
$M\times (0,T)$. 

Let $\Omega\subset M$ be a bounded domain with smooth boundary $\1\Omega$.
Then there exists a constant $H>0$ such that the second fundamental form II 
of $\1\Omega$ with respect to the unit outward normal $\1/\1\nu$ of $\1\Omega$
and metric $g(0)$ is uniformly bounded below by $-H$. 
For any $x\in\Omega$ let $\rho^t(x)$ be the distance of $x$ from $\1\Omega$ 
with respect to the metric $g(t)$ and $\rho (x)=\rho^0(x)$. Note that by
(0.1) and (4.1) there exist constants $c_1>0$, $c_2>0$, such that
$$\left\{\aligned
&c_1g_{ij}(x,t_1)\le g_{ij}(x,t_2)\le c_2g_{ij}(x,t_1)\quad\forall x\in M,
0\le t_1,t_2\le T\\
&c_1g^{ij}(x,t_1)\le g^{ij}(x,t_2)\le c_2g^{ij}(x,t_1)\quad\forall x\in M,
0\le t_1,t_2\le T\\
&c_1\rho^{t_1}(x)\le\rho^{t_2}(x)\le c_2\rho^{t_1}(x)\qquad\qquad
\forall x\in\Omega,0\le t_1,t_2\le T\\
&c_1dV_{t_1}\le dV_{t_2}\le c_2dV_{t_1}\qquad\qquad\qquad\text{ in }M\quad
\forall 0\le t_1,t_2\le T.\endaligned\right.\tag 4.2
$$
For any $\delta>0$, let $\Omega_{\delta}=\{x\in\Omega:\rho (x)\ge\delta\}$.
Let $p\in\Omega$ and $\2{u}(x,t)$ be the Dirichlet fundamental solution 
of the conjugate heat equation (0.3) in $\Omega\times (0,T)$ which 
satisfies (0.4) with
$$
\2{u}=0\quad\text{ on }\1\Omega\times (0,T).\tag 4.3
$$ 
Let $\2{f}$, $\2{v}$, be given by
$$\left\{\aligned
&\2{u}(x,t)=\frac{e^{-\2{f}(x,t)}}{(4\pi\tau)^{\frac{n}{2}}}\\
&\2{v}=[\tau (2\Delta^t\2{f}-|\nabla^t\2{f}|^2+R)+\2{f}-n]\2{u}
\endaligned\right.
$$
and let $\4{u}(x,\tau)=\2{u}(x,T-t)$, $d\4{V}_{\tau}=dV_t$, where $\tau=T-t$.
We choose $0<\delta<1$ such that $p\in\Omega_{3\delta}$ (cf. \cite{C}),
$$\left\{\aligned
&\sqrt{k_0}\tan (3\delta\sqrt{k_0})\le\frac{H}{2}+\frac{1}{2}\\
&\frac{H}{\sqrt{k_0}}\tan (3\delta\sqrt{k_0})\le\frac{1}{2}.
\endaligned\right.\tag 4.4
$$
By the maximum principle,
$$
\2{u}(x,t)\le\Cal{Z}(p,T,x,t)\quad\text{ in }\Omega\times (0,T).\tag 4.5
$$

By compactness and an argument similar to the proof of Corollary 4.1 of 
\cite{CTY} and Theorem 1.2 we have

\proclaim{\bf Theorem 4.1}
For any $\alpha>1$, $\3>0$, and $0<\delta_2\le\delta$, there exists a 
constant $C_1>0$ depending on $k_0$, $\alpha$, $\3$ and $\delta_2$ such that
$$
\frac{|\nabla^t\4{u}|^2}{\4{u}^2}-\alpha\frac{\4{u}_{\tau}}{\4{u}}\le C_1
+\frac{n(1+\3)}{2\tau}\quad\text{ in }\2{\Omega}_{\delta_2}\times (0,T]
\tag 4.6
$$
where $t=T-\tau$.
\endproclaim

\proclaim{\bf Lemma 4.2}
There exist a constant $0<\tau_0<\delta^2$ and  constants $C_2>0$, $C_3>0$,
and $D>1$ independent of $\tau_0$ such that for any $0<\tau_1\le\tau_0$,
$$\aligned
(i)&\quad 0<\4{u}(x,\tau)\le\frac{C_2}{\tau_1^\frac{n}{2}}
e^{-\frac{\delta^2}{D\tau_1}}\quad\forall x\in\Omega\setminus
\Omega_{2\delta}, 0<\tau\le\tau_1\\
(ii)&\quad 0<\4{u}(x,\tau)\le C_3
\quad\forall x\in\Omega\setminus\Omega_{2\delta},0<\tau\le T\\
(iii)&\quad|\nabla^t\4{u}(x,\tau)|+|\nabla^t\nabla^t\4{u}(x,\tau)|
\le C_3\quad\forall x\in\2{\Omega}\setminus\Omega_{2\delta},0<\tau\le T
\endaligned
$$
where $t=T-\tau$.
\endproclaim
\demo{Proof}
The left hand side of (i) and (ii) follows by the strong maximum priniciple.
By Corollary 5.2 of \cite{CTY} there exist constants $C>0$ and $D>1$
such that
$$
\Cal{Z}(p,T;x,t)\le\frac{C}{V_p(\sqrt{\tau})}e^{-\frac{r^2(p,x)}{D(T-t)}}.
\quad\forall x\in M,0\le t<T, \tau=T-t.\tag 4.7
$$
By the same argument as the proof of Lemma 7.6 of \cite{CTY} there exist 
constants $C_1>0$ and $C_2>0$ such that
$$
C_1\tau^{\frac{n}{2}}\le V_p^t(\sqrt{\tau})\le C_2\tau^{\frac{n}{2}}\quad
\forall 0<\tau\le T, t=T-\tau.\tag 4.8
$$
Hence by (4.2), (4.5), (4.7) and (4.8),
$$\align
&\4{u}(x,\tau)\le\frac{C}{V_p(\sqrt{\tau})}e^{-\frac{r^2(p,x)}{D\tau}}
\le\frac{C'}{V_p^t(\sqrt{\tau})}e^{-\frac{r^2(p,x)}{D\tau}}
\le\frac{C_2}{\tau^{\frac{n}{2}}}e^{-\frac{r^2(p,x)}{D\tau}}
\quad\forall x\in\2{\Omega},0<\tau\le T\\
\Rightarrow\quad&\4{u}(x,\tau)\le\frac{C_2}{\tau^{\frac{n}{2}}}
e^{-\frac{\delta^2}{D\tau}}\quad\forall x\in\2{\Omega}\setminus
\Omega_{2\delta},0<\tau\le T\tag 4.9
\endalign
$$
for some constant $C_2>0$. Let $H(\tau)=\tau^{-\frac{n}{2}}
e^{-\frac{\delta^2}{D\tau}}$, $\tau_0=\delta^2/(nD)$ and
$0<\tau_1\le\tau_0$. Then $H'(\tau)\ge 0$ for all $0<\tau\le\tau_0$. Hence 
by (4.9) (i) follows. By (4.9) and (i) we get (ii).

We now extend $\4{u}$ to a function on $(\2{\Omega}\setminus\Omega_{2\delta})
\times (-\infty,T]$ by setting $\4{u}=0$ on $(\2{\Omega}\setminus
\Omega_{2\delta})\times (-\infty,0)$. We also extend $g_{ij}$ to a metric on
$(\2{\Omega}\setminus\Omega_{2\delta})\times [0,\infty)$ by setting 
$g_{ij}(x,t)=g_{ij}(x,T)$ for all $t\ge T$. Then $\4{u}$ is a non-negative
solution of 
$$
\4{u}_{\tau}=\Delta^t\4{u}-R(x,t)\4{u}\quad\text{ in }
\2{\Omega}\setminus\Omega_{2\delta}\times (-\infty,T],t=T-\tau.\tag 4.10
$$
Hence by the parabolic regularity theory \cite{LSU} $\4{u}\in 
C^{\infty}(\2{\Omega}\setminus\Omega_{2\delta}\times (-\infty,T])$.
Thus
$$\left\{\aligned
&|\nabla^t\4{u}(x,\tau)|\le\max\Sb y\in\2{\Omega}\setminus\Omega_{2\delta}\\
0<\tau\le T\endSb |\nabla^t\4{u}(y,\tau)|<\infty
\quad\forall x\in\2{\Omega}\setminus\Omega_{2\delta},0<\tau\le T\\
&|\nabla^t\nabla^t\4{u}(x,\tau)|\le\max\Sb y\in\2{\Omega}\setminus
\Omega_{2\delta}\\0<\tau\le T\endSb |\nabla^t\nabla^t\4{u}(y,\tau)|<\infty
\quad\forall x\in\2{\Omega}\setminus\Omega_{2\delta},0<\tau\le T
\endaligned\right.
$$
and (iii) follows.
\enddemo

\proclaim{\bf Theorem 4.3}
Let $0<\tau_0<\delta^2$ and $C_2>0$ be as given in Lemma 4.2. 
Then there exists a constant $C_4>0$ depending on $k_0$ such that
$$\aligned
(i)&\quad\frac{|\nabla^t\4{u}(x,\tau)|}{\4{u}(x,\tau)}
\le\frac{C_4}{\rho (x)}\biggl (1+\log\biggl (C_2\frac{e^{-\frac{\delta^2}
{D\tau}}}{\tau^{n/2}\4{u}}
\biggr )\biggr )\quad\forall x\in\Omega,\rho (x)<\sqrt{\tau}\\
(ii)&\quad\frac{|\nabla^t\4{u}(x,\tau)|}{\4{u}(x,\tau)}
\le\frac{C_4}{\sqrt{\tau}}\biggl (1+\log\biggl (C_2
\frac{e^{-\frac{\delta^2}{D\tau}}}{\tau^{n/2}\4{u}}
\biggr )\biggr )\quad\forall x\in\Omega,\sqrt{\tau}\le\rho (x)\le\delta
\endaligned
$$
holds for any $0<\tau\le\tau_0$ where $t=T-\tau$ and $D>1$ is as given in 
Lemma 4.2.
\endproclaim
\demo{Proof}
Let $0<\tau\le\tau_0$. We divide the proof into two cases.

\noindent $\underline{\text{Case 1}}$: $\rho (x)<\sqrt{\tau}$

By applying Theorem 3.1 to the domain $Q_1=Q_{\rho (x),\frac{\rho (x)^2}{4}}
(x,\tau)$ there exists a constant $C_1>0$ such that
$$
\frac{|\nabla^t\4{u}(x,\tau)|}{\4{u}(x,\tau)}
\le C_1\biggl (1+\frac{1}{\rho (x)}\biggr )\biggl (
1+\log\biggl (\frac{A_1}{\4{u}}\biggr )\biggr )\tag 4.11
$$
holds where $A_1=\sup_{Q_1}\4{u}$.

\noindent $\underline{\text{Case 2}}$: $\sqrt{\tau}\le\rho (x)\le\delta$

By applying Theorem 3.1 to the domain $Q_2=Q_{\sqrt{\frac{\tau}{2}},
\frac{\tau}{2}}(x,\tau)$ there exists a constant $C_1>0$ such that
$$
\frac{|\nabla^t\4{u}(x,\tau)|}{\4{u}(x,\tau)}
\le C_1\biggl (1+\frac{1}{\sqrt{\tau}}\biggr )\biggl (
1+\log\biggl (\frac{A_2}{\4{u}}\biggr )\biggr )\tag 4.12
$$
holds where $A_2=\sup_{Q_2}\4{u}$.
By Lemma 4.2 there exist constants $C_2>0$, $D>1$, such that
$$
A_1,A_2\le\frac{C_2}{\tau^\frac{n}{2}}e^{-\frac{\delta^2}{D\tau}}\tag 4.13
$$
holds for any $0<\tau\le\tau_0$. Hence by (4.11), (4.12) and (4.13), the 
lemma follows.
\enddemo

By a similar argument we have

\proclaim{\bf Theorem 4.4}
Let $C_3>0$ be as given in Lemma 4.2. Then there exists a constant $C_4>0$ 
depending on $k_0$ such that
$$\aligned
(i)&\quad\frac{|\nabla^t\4{u}(x,\tau)|}{\4{u}(x,\tau)}
\le\frac{C_4}{\rho (x)}\biggl (1+\log\biggl (\frac{C_3}{\4{u}}
\biggr )\biggr )\quad\forall x\in\Omega,\rho (x)<\sqrt{\tau}\\
(ii)&\quad\frac{|\nabla^t\4{u}(x,\tau)|}{\4{u}(x,\tau)}
\le\frac{C_4}{\sqrt{\tau}}\biggl (1+\log\biggl (\frac{C_3}{\4{u}}
\biggr )\biggr )\quad\forall x\in\Omega,\sqrt{\tau}\le\rho (x)\le\delta
\endaligned
$$
holds for any $0<\tau\le T$ where $t=T-\tau$.
\endproclaim

By Lemma 4.2 and Theorem 4.3 we have the following corollary.

\proclaim{\bf Corollary 4.5}
Let $0<\tau_0<\delta^2$ and $D>1$ be as given in Lemma 4.2. Then for any
$a>0$ there exists a constant $C>0$ depending on $k_0$ and $a$ such that
$$\aligned
(i)&\quad\frac{|\nabla^t\4{u}(x,\tau)|}{\4{u}(x,\tau)}
\le\frac{C}{\rho (x)}\biggl (1+\biggl (
\frac{e^{-\frac{\delta^2}{D\tau}}}{\tau^{\frac{n}{2}}\4{u}}\biggr )^a\biggr )
\quad\forall x\in\Omega,\rho (x)<\sqrt{\tau}\\
(ii)&\quad\frac{|\nabla^t\4{u}(x,\tau)|}{\4{u}(x,\tau)}
\le\frac{C}{\sqrt{\tau}}\biggl (1+\biggl (
\frac{e^{-\frac{\delta^2}{D\tau}}}{\tau^{\frac{n}{2}}\4{u}}\biggr )^a
\biggr )\quad\forall x\in\Omega,\sqrt{\tau}\le\rho (x)\le\delta
\endaligned
$$
holds for any $0<\tau\le\tau_0$.
\endproclaim

\proclaim{\bf Lemma 4.6}
Let $0<\tau_0<\delta^2$ and $D>1$ be as given in Lemma 4.2. Then there 
exists a constant $C_5>0$ depending on $k_0$ such that
$$\aligned
(i)&\quad\frac{|\nabla^t\4{u}(x,\tau)|^2}{\4{u}(x,\tau)}
\le C_5\frac{e^{-\frac{\delta^2}{4D\tau}}}
{\tau^{\frac{3n}{8}}\rho (x)^{\frac{7}{4}}}
\quad\forall x\in\Omega,\rho (x)<\sqrt{\tau}\\
(ii)&\quad\frac{|\nabla^t\4{u}(x,\tau)|^2}{\4{u}(x,\tau)}
\le C_5\frac{e^{-\frac{\delta^2}{2D\tau}}}{\tau^{\frac{n}{2}+1}}
\quad\forall x\in\Omega,\sqrt{\tau}\le\rho (x)\le\delta
\endaligned
$$
holds for any $0<\tau\le\tau_0$.
\endproclaim
\demo{Proof}
By Corollary 4.5 there exists a constant $C>0$ such that
$$
\frac{|\nabla^t\4{u}(x,\tau)|^2}{\4{u}(x,\tau)}
\le\left\{\aligned
&\frac{C}{\rho (x)}\biggl (1+\frac{1}{(\tau^{n/2}\4{u})^{1/4}}\biggr )
|\nabla\4{u}|\quad\forall\rho (x)<\sqrt{\tau}\\
&\frac{C}{\sqrt{\tau}}\biggl (1+\frac{1}{(\tau^{n/2}\4{u})^{1/4}}\biggr )
|\nabla\4{u}|
\quad\forall\sqrt{\tau}\le\rho (x)\le\delta\endaligned\right.\tag 4.14
$$
holds for any $0<\tau\le\tau_0$ where $t=T-\tau$.
We will now let $C>0$ be a generic constant that may change from line to line.
By Lemma 4.2 and Corollary 4.5,
$$\align
\frac{|\nabla^t\4{u}|}{\4{u}^{\frac{1}{4}}}
\le&C\frac{\4{u}^{3/4}}{\sqrt{\tau}}
\biggl (1+\frac{1}{(\tau^{n/2}\4{u})^{1/4}}\biggr )
=C\biggl (\frac{\4{u}^{3/4}}{\sqrt{\tau}}
+\frac{\4{u}^{1/2}}{\tau^{\frac{n}{8}+\frac{1}{2}}}\biggr )\\
\le&C\biggl (\frac{1}{\sqrt{\tau}}\biggl (
\frac{e^{-\frac{\delta^2}{D\tau}}}{\tau^{n/2}}
\biggr )^{3/4}+\frac{1}{\tau^{\frac{n}{8}+\frac{1}{2}}}\biggl (
\frac{e^{-\frac{\delta^2}{D\tau}}}{\tau^{n/2}}\biggr )^{1/2}\biggr )\\
\le&\frac{C}{\tau^{\frac{3n}{8}+\frac{1}{2}}}e^{-\frac{\delta^2}{2D\tau}}
\qquad\qquad\qquad\qquad\forall x\in\Omega_{\sqrt{\tau}}\setminus
\Omega_{\delta},0<\tau\le\tau_0\tag 4.15
\endalign
$$
and
$$\align
\frac{|\nabla^t\4{u}|}{\4{u}^{\frac{1}{4}}}
\le&C\frac{\4{u}^{3/4}}{\rho (x)}
\biggl (1+\frac{1}{(\tau^{n/2}\4{u})^{1/4}}\biggr )
=\frac{C}{\rho (x)^{3/4}}\biggl (\frac{\4{u}}{\rho (x)}\biggr )^{1/4}
\biggl (\4{u}^{1/2}
+\frac{\4{u}^{1/4}}{\tau^{\frac{n}{8}}}\biggr )\\
\le&\frac{C}{\rho (x)^{3/4}}\biggl (\frac{\4{u}}{\rho (x)}\biggr )^{1/4}
\biggl (\biggl (\frac{e^{-\frac{\delta^2}{D\tau}}}{\tau^{n/2}}
\biggr )^{1/2}+\frac{1}{\tau^{\frac{n}{8}}}\biggl (
\frac{e^{-\frac{\delta^2}{D\tau}}}{\tau^{n/2}}\biggr )^{1/4}\biggr )\\
\le&\frac{C}{\rho (x)^{3/4}}\biggl (\frac{\4{u}}{\rho (x)}\biggr )^{1/4}
\cdot\frac{e^{-\frac{\delta^2}{4D\tau}}}{\tau^{n/4}}\qquad
\qquad\forall x\in\Omega\setminus\Omega_{\sqrt{\tau}},0<\tau\le\tau_0.
\tag 4.16\endalign
$$
By (4.14), (4.15), and Lemma 4.2,
$$\align
\frac{|\nabla^t\4{u}(x,\tau)|^2}{\4{u}(x,\tau)}
\le&\frac{C}{\sqrt{\tau}}\biggl (\frac{1}{\tau^{\frac{3n}{8}+\frac{1}{2}}}
e^{-\frac{\delta^2}{2D\tau}}\4{u}^{\frac{1}{4}}
+\frac{1}{\tau^{\frac{n}{2}+\frac{1}{2}}}e^{-\frac{\delta^2}{2D\tau}}\biggr )
\quad\forall\sqrt{\tau}\le\rho (x)\le\delta,0<\tau\le\tau_0\\
\le&\frac{C}{\sqrt{\tau}}\biggl (\frac{1}{\tau^{\frac{3n}{8}+\frac{1}{2}}}
+\frac{1}{\tau^{\frac{n}{2}+\frac{1}{2}}}\biggr )
e^{-\frac{\delta^2}{2D\tau}}\qquad\qquad\quad\forall\sqrt{\tau}\le\rho (x)
\le\delta,0<\tau\le\tau_0\\
\le&C\frac{e^{-\frac{\delta^2}{2D\tau}}}{\tau^{\frac{n}{2}+1}}
\qquad\qquad\qquad\qquad\qquad\qquad\qquad\forall\sqrt{\tau}\le\rho (x)
\le\delta,0<\tau\le\tau_0
\endalign
$$ 
where $t=T-\tau$ and (ii) follows. 

By (4.4) (cf. \cite{Wa}, \cite{C}, and \cite{Ch}) 
for any $x\in\Omega\setminus\Omega_{\delta}$, there exists a unique 
normalized minimizing geodesic $\gamma_x:[0,\rho (x)]\to\2{\Omega}$ with 
respect to the metric $g(0)$ such that $\gamma_x(0)\in\1\Omega$, 
$\gamma_x(\rho (x))=x$, and $\gamma_x'(0)$ is perpendicular to the tangent 
plane $T_{\gamma_x(0)}(\1\Omega)$ at $\gamma_x(0)$. 
By (4.2) and Lemma 4.2 for any $x\in\Omega\setminus\Omega_{\delta}$, 
$0<\tau\le\tau_0$,
$$
\frac{\4{u}(x,\tau)}{\rho (x)}=\frac{\int_0^{\rho (x)}
\frac{\1}{\1 s}\4{u}(\gamma_x(s),\tau)\,ds}{\rho (x)}
\le C\sup\Sb y\in\Omega\\0<\tau\le\tau_0\endSb|\nabla^0\4{u}(y,\tau)|
\le C\sup\Sb y\in\Omega\\0<\tau\le\tau_0\endSb|\nabla^t\4{u}(y,\tau)|
\le C<\infty.\tag 4.17
$$
By (4.14), (4.16), (4.17) and Lemma 4.2 we get (i) and the lemma follows.
\enddemo

By Lemma 4.2, Theorem 4.4, and an argument similar to the proof of Lemma 4.6 
we have

\proclaim{\bf Theorem 4.7}
Then there exists a constant $C_6>0$ depending on $k_0$ such that
$$\align
(i)&\quad\frac{|\nabla^t\4{u}(x,\tau)|^2}{\4{u}(x,\tau)}
\le\frac{C_6}{\rho (x)}\biggl (1+\frac{1}{\rho (x)^{\frac{3}{4}}}\biggr )
\quad\forall x\in\Omega,\rho (x)<\sqrt{\tau}\\
(ii)&\quad\frac{|\nabla^t\4{u}(x,\tau)|^2}{\4{u}(x,\tau)}
\le\frac{C_6}{\sqrt{\tau}}\biggl (1+\frac{1}{\sqrt{\tau}}
\biggr )\quad\forall x\in\Omega,\sqrt{\tau}\le\rho (x)\le\delta
\endalign
$$
holds for any $0<\tau\le T$.
\endproclaim


\Refs

\ref
\key CTY\by\ A.~Chau, L.F.~Tam and C.~Yu\paper Pseudolocality 
for the Ricci flow and applications,\ \ \linebreak
http://arxiv.org/abs/math/0701153\endref 

\ref
\key Ch\by I.~Chavel\book Riemannian geometry: A modern introduction
\publ Cambridge University Press\publaddr Cambridge, United Kingdom
\yr 1995\endref

\ref
\key C\by R.~Chen\paper Neumann eigenvalue estimate on a compact Riemannian 
manifold\jour Proc. AMS\vol 108\yr 1990\pages 961--970\endref

\ref
\key CLN\by \ \ \ B.~Chow and P.~Lu and Lei Ni\book Hamilton's Ricci flow,
Graduate Studies in Mathematics, vol. 77\publ Amer. Math. Soc.
\publaddr Providence, R.I., U.S.A.\yr 2006\endref

\ref 
\key H\by R.S.~Hamilton\paper The formation of singularities in the Ricci 
flow\jour Surveys in differential geometry, Vol. II (Cambridge, MA, 1993),
7--136, International Press, Cambridge, MA, 1995\endref

\ref 
\key Hs1\by \ \ S.Y.~Hsu, Uniqueness of solutions of Ricci flow
on complete noncompact manifolds, http://arxiv.org\linebreak
/abs/0704.3468\endref

\ref
\key Hs2\by \ S.Y.~Hsu\paper Maximum principle and convergence of
fundamental solutions for the Ricci flow,\linebreak 
http://arxiv.org/abs/math/0711.1236\endref

\ref
\key KZ\by \ S.~Kuang and Q.S.~Zhang\paper A gradient estimate for
all positive solutions of the conjugate heat equation under Ricci
flow, http://arxiv.org/abs/math/0611298\endref

\ref
\key LSU\by \ \ O.A.~Ladyzenskaya, V.A.~Solonnikov, and
N.N.~Uraltceva\book Linear and quasilinear equations of
parabolic type\publ Transl. Math. Mono. Vol 23,
Amer. Math. Soc.\publaddr Providence, R.I.\yr 1968\endref

\ref
\key LY\by P.~Li and S.T.~Yau\paper On the parabolic kernel of the 
Schr\"odinger operator\jour Acta Math.\vol 156\yr 1986\pages 153--201
\endref
 
\ref
\key N1\by L.~Ni\paper The entropy formula for linear heat equation
\jour J. Geometric Analysis\vol 14(1)\yr 2004\pages 87--100
\endref

\ref
\key N2\by L.~Ni\paper Addenda to ``The entropy formula for linear heat 
equation''\jour J. Geometric Analysis\vol 14(2)\yr 2004\pages 369--374
\endref

\ref
\key N3\by L.~Ni\paper A note on Perelman's LYH-type inequality\jour Comm.
Anal. and Geom.\vol 14(5)\yr 2006\pages 883-905\endref

\ref
\key P\by G.~Perelman\paper The entropy formula for the Ricci flow and its 
geometric applications,\linebreak http://arXiv.org/abs/math.DG/0211159\endref 

\ref
\key SY\by R.~Schoen and S.T.~Yau\book Lectures on differential geometry,
in 'Conference proceedings and Lecture Notes in Geometryand Topology', 1
\publ International Press\publaddr \yr 1994\endref

\ref
\key SZ\by P.~Souplet and Q.S.~Zhang\paper Sharp gradient estimate and
Yau's Liouville theorem for the heat equation on noncompact manifolds
\jour Bull. London math. Soc.\vol 38\yr 2006\pages 1045--1053\endref

\ref
\key W\by J.~Wang\paper Global heat kernel estimates\jour Pacific J. Math.
\vol 178(2)\yr 1997\pages 377--398
\endref

\ref
\key Wa\by \ F.W.~Warner\paper Extension of the Rauch comparison theorem to
submanifolds\jour Trans. Amer. Math. Soc.\vol 122\yr 1966\pages 341--356
\endref

\ref
\key Z\by Q.S.~Zhang\paper Some gradient estimates for the heat equation
on domains and for an equation by Perelman\jour Int. Math. Res. Notice
\yr 2006\pages Art. ID 92314, 39 pp\endref


\endRefs
\enddocument